\documentclass{amsart}
\usepackage{amsthm}
\usepackage{amssymb}
\usepackage{upref}
\usepackage{amscd}
\usepackage{xspace}

\numberwithin{equation}{section}
\newtheorem{main}{Theorem}

\newtheorem{theorem}[equation]{Theorem}
\newtheorem{lemma}[equation]{Lemma}
\newtheorem{proposition}[equation]{Proposition}
\newtheorem{corollary}[equation]{Corollary}

\theoremstyle{definition}
\newtheorem{definition}[equation]{Definition}
\newtheorem*{remark}{Remark}

\DeclareMathOperator{\Ann}{Ann}

\DeclareMathOperator{\Hom}{Hom}

\DeclareMathOperator{\Ext}{Ext}

\DeclareMathOperator{\colim}{colim}

\renewcommand{\smash}{\wedge}

\newcommand{\Cech}{\check{C}}
\newcommand{\CechH}{\check{C}H}
\newcommand{\textCech}{\v{C}ech\xspace}
\newcommand{\bt}{\bullet}

\newcommand{\cat}[1]{\mathcal{#1}}

\newcommand{\Z}{\mathbb{Z}}

\newcommand{\comod}{\text{-comod}}

\newcommand{\mathcolon}{\colon\,}

\newcommand{\ulp}{\textup{(}}
\newcommand{\urp}{\textup{)}}

\begin{document}
 
\title{Local cohomology of $BP_{*}BP$-comodules} 

\date{\today}

\author{Mark Hovey}
\address{Department of Mathematics \\ Wesleyan University
\\ Middletown, CT 06459}
\email{hovey@member.ams.org}

\author{Neil Strickland}
\address{Department of Pure Mathematics \\ University of Sheffield
\\ Sheffield S3 7RH \\
England}
\email{N.P.Strickland@sheffield.ac.uk}


\begin{abstract}
Given a spectrum $X$, we construct a spectral sequence of
$BP_{*}BP$-comodules that converges to $BP_{*}(L_{n}X)$, where $L_{n}X$
is the Bousfield localization of $X$ with respect to the Johnson-Wilson
theory $E(n)_{*}$.  The $E_{2}$-term of this spectral sequence consists
of the derived functors of an algebraic version of $L_{n}$.  We show how
to calculate these derived functors, which are closely related to local
cohomology of $BP_{*}$-modules with respect to the ideal $I_{n+1}$.  
\end{abstract}

\maketitle

\section*{Introduction}\label{sec-intro}

The most common approach to understanding stable homotopy theory
involves localization.  One first localizes at a fixed prime $p$; after
doing so there is a tower of localization functors 
\[
\dotsb L_{n} \xrightarrow{} L_{n-1} \xrightarrow{} \dotsb \xrightarrow{}
L_{1}\xrightarrow{} L_{0}.  
\]
Each functor $L_{n}$ retains a little more information than the previous
one $L_{n-1}$; the homotopy inverse limit of the $L_{n}X$ is $X$ itself
for a finite $p$-local spectrum $X$.  These localization functors come
from the Brown-Peterson homology theory $BP$, where
\[
BP_{*}(S^{0})\cong \Z_{(p)}[v_{1},v_{2},\dotsc]
\]
with $|v_{i}|=2(p^{i}-1)$.  The generators $v_{n}$ are not well-defined,
but the ideals $I_{n}=(p,v_{1},\dotsc ,v_{n-1})$ for $0\leq n\leq
\infty$ are well-defined.  The functor $L_{n}$ is Bousfield
localization with respect to $v_{n}^{-1}BP$; all the different choices
for $v_{n}$ give the same localization functor.  

In previous work~\cite{hovey-strickland-comodules}, the authors
constructed an algebraic endofunctor $L_{n}$ on the category of
$BP_{*}BP$-comodules, analogous to the chromatic localization $L_{n}$ on
spectra.  This functor $L_{n}$ is the localization obtained by inverting
all maps of comodules whose kernel and cokernel are $v_{n}$-torsion (or,
equivalently, $I_{n+1}$-torsion).  The $L_{n}$-local comodules are
equivalent to the category of $E(n)_{*}E(n)$-comodules, or to the
category of $E_{*}E$-comodules for any Landweber exact commutative ring
spectrum with $E_{*}/I_{n+1}=0$ but $E_{*}/I_{n}\neq 0$.

In~\cite{hovey-strickland-comodules}, our main interest was algebraic.
In this paper, we compare our algebraic version of $L_{n}$ with the
topological one.  As always, when one has a topological version of an
algebraic construction, one expects a spectral sequence converging to
the topological construction whose $E_{2}$-term involves the derived
functors of the algebraic construction.  Since the algebraic $L_{n}$ is
left exact, it has right derived functors $L_{n}^{i}$.  We prove the
following theorem.

\begin{main}\label{main-A}
Let $X$ be a spectrum.  There is a natural spectral sequence
$E_{*}^{**}(X)$ with $d_{r}\mathcolon
E_{r}^{s,t}\xrightarrow{}E_{r}^{s+t,t+r-1}$ and $E_{2}$-term
$E_{2}^{s,t}(X)\cong (L^{s}_{n}BP_{*}X)_{t}$, converging conditionally
and strongly to $BP_{t-s}(L_{n}X)$.  This is a spectral sequence of
$BP_{*}BP$-comodules, in the sense that $E_{r}^{s,*}$ is a graded
$BP_{*}BP$-comodule for all $r\geq 2$ and $d_{r}\mathcolon
E_{r}^{s,*}\xrightarrow{}E_{r}^{s+r,*}$ is a $BP_{*}BP$-comodule map of
degree $r-1$.  Furthermore, every element in $E_{2}^{0,*}$ that comes
from $BP_{*}X$ is a permanent cycle.
\end{main}

For this spectral sequence to be useful, we need to be able to compute
the $E_{2}$-term.  The derived functors $L_{n}^{i}$ turn out to be
closely related to local cohomology, which is well-known in commutative
algebra and was introduced to algebraic topology by
Greenlees~\cite{greenlees-spectral}.  Recall
from~\cite{greenlees-may-completions} that, given an ideal $I$ in a ring
$R$, one can form the local cohomology $H_{I}^{*}(-)$ and the Cech
cohomology $\CechH_{I}^{*}(-)$ of an $R$-module $M$.  Although it is
not phrased this way in~\cite{greenlees-may-completions}, the functor
$\CechH_{I_{n+1}}^{0}$ on the category of $BP_{*}$-modules is the
localization functor that inverts all maps of modules whose kernel and
cokernel are $v_{n}$-torsion.  Thus $\CechH_{I_{n+1}}^{0}$ is the
analog of $L_{n}$ in the category of $BP_{*}$-modules, and therefore
Cech cohomology is simply the derived functors of localization on the
category of $BP_{*}$-modules.

\begin{main}\label{main-B}
Let $M$ be a $BP_{*}BP$-comodule. 
\begin{enumerate}
\item [(1)] We have $L_{n}^{i}M\cong \CechH_{I_{n+1}}^{i}M$. 
\item [(2)] $L_{n}^{i}(M)=0$ for $i>n$.  
\item [(3)] $L_{n}^{i}(M)$ is $I_{n+1}$-torsion for all $i>0$.  
\item [(4)] If $v_{j}$ acts isomorphically on $M$ for some $0\leq j\leq
n$, then $L_{n}M=M$ and $L_{n}^{i}M=0$ for $i>0$.  
\item [(5)] If $k<n$, then $L_{n}^{i}(BP_{*}/I_{k})=0$ unless $i=0$ or
$i=n-k$.  We have 
\[
L_{n}(BP_{*}/I_{k}) = BP_{*}/I_{k}
\]
and 
\[
L_{n}^{n-k}(BP_{*}/I_{k}) = BP_{*}/(p, v_{1},\dotsc
,v_{k-1},v_{k}^{\infty},\dotsc ,v_{n}^{\infty}).  
\]
\item [(6)] We have $L_{n}^{i}(BP_{*}/I_{n})=0$ for $i>0$ and 
\[
L_{0}(BP_{*}) = BP_{*}\otimes \mathbb{Q} \text{ and }
L_{n}(BP_{*}/I_{n})=v_{n}^{-1}BP_{*}/I_{n} \text{ for } n>0.
\]
\item [(7)] If $k>n$, then $L_{n}^{i}(BP_{*}/I_{k})=0$ for all $i$.  
\item [(8)] $L_{n}^{i}$ commutes with filtered colimits, arbitrary
direct sums, and finite limits.  
\end{enumerate}
\end{main}

Most of Theorem~\ref{main-B}, except part~(8), would follow from part~(1) of
it and known facts about local cohomology.  However, local cohomology is
generally considered only for Noetherian rings, and $BP_{*}$ is not
Noetherian.  This turns out not to be a problem, but because there is no
discussion of non-Noetherian local cohomology in the literature, and
because it is not very hard, we offer direct proofs of the remaining
parts of Theorem~\ref{main-B}.  

In the light of Theorem~\ref{main-B} the reader may naturally wonder
whether there is a connection between the local cohomology spectral
sequence of~\cite{greenlees-spectral}
and~\cite{greenlees-may-completions} and our spectral sequence.
Recall that Greenlees and May begin with the category of modules over
a strictly commutative ring spectrum; since $BP$ is not known to be
such, we must begin with $MU$.  Combining Theorems~5.1 and~6.1
of~\cite{greenlees-may-completions}, and applying them to the
$MU$-module spectrum $MU\smash X$, then gives a spectral sequence
converging to $MU_{*}(L_{n}X)$ whose $E_{2}$ is
$\CechH_{I_{n+1}}^{-s,-t}(MU_{*}X)$.  Our spectral sequence would
coincide with this Greenlees-May spectral sequence if we used $MU$
instead of $BP$ (and we reindexed the spectral sequence).  However,
our construction allows us to conclude that we have a spectral
sequence \textbf{of comodules}, which the Greenlees-May construction
does not.  This significantly restricts the possible differentials and
extensions that can occur in the spectral sequence.

The authors would like to thank the Isaac Newton Institute for
Mathematical Sciences for its support during our collaboration.  We
also thank John Greenlees for many helpful conversations on the subject
matter of this paper.  

\section{The functor $L_{n}$}\label{sec-ln}

In this section, we define our localization functor $L_{n}$ and derive
some basic properties of it.  The objective of this section, then, is to
prove those parts of Theorem~\ref{main-B} dealing with the localization
functor $L_{n}$ itself, but not its derived functors.  We do not quite
succeed in this, as Theorem~\ref{main-B}(1) requires some different ideas
that we postpone to Section~\ref{sec-Cech}.

In this section, and throughout the paper, $n$ will be a fixed
nonnegative integer, and
\[
\Phi \mathcolon (BP_{*},BP_{*}BP)\xrightarrow{}(E(n)_{*}, E(n)_{*}E(n))
\]
will denote the evident map of Hopf algebroids.  There is an induced
exact functor 
\[
\Phi_{*}\mathcolon BP_{*}BP\comod \xrightarrow{} E(n)_{*}E(n)\comod 
\]
of the categories of graded comodules that takes $M$ to
$E(n)_{*}\otimes_{BP_{*}}M$.  As explained
in~\cite[Proposition~1.2.3]{hovey-comodule-homotopy}, $\Phi_{*}$ has a
left exact right adjoint
\[
\Phi^{*} \mathcolon E(n)_{*}E(n)\comod \xrightarrow{} BP_{*}BP\comod.
\]
We prove in Section~2 of~\cite{hovey-strickland-comodules} that
$\Phi^{*}$ is a fully faithful embedding with 
\[
\Phi_{*}\Phi^{*}M\cong M,
\]
and in Section~4 of that paper that the composite functor
$\Phi^{*}\Phi_{*}$ is localization with respect to the hereditary
torsion theory consisting of $v_{n}$-torsion comodules.  It is this
composite $\Phi^{*}\Phi_{*}$ that we denote by $L_{n}$.  For a quick
review of the theory of localization with respect to hereditary torsion
theories, see the discussion immediately following
Corollary~\ref{cor-Cech-local}.  Since the $v_{n}$-torsion comodules are
the smallest hereditary torsion theory containing
$BP_{*}/I_{n+1}$~\cite{hovey-strickland-comodules}, $L_{n}$ is
localization away from $BP_{*}/I_{n+1}$ so is analogous to $L_{n}^{f}$
on the category of spectra.  On the other hand, the $v_{n}$-torsion
comodules are precisely the kernel of
$\Phi_{*}$~\cite{hovey-strickland-comodules}, so $L_{n}$ is also
analogous to the functor $L_{n}$ on the category of spectra.

Because the collection of $v_{n}$-torsion comodules is a hereditary
torsion theory, the submodule $T_{n}M$ of $v_{n}$-torsion elements in a
comodule $M$ is in fact a subcomodule.  This also follows
from~\cite[Corollary~2.4]{johnson-yosimura} 
or~\cite[Corollary 2]{landweber-comodules}. 

The most basic facts about $L_{n}$ are contained in the following
proposition.  

\begin{proposition}\label{prop-Ln-basic}
\begin{enumerate}
\item [(a)] $L_{n}$ is left exact, idempotent, and commutes with finite
limits.
\item [(b)] Given a map $f$ of comodules, $L_{n}f$ is an isomorphism if
and only if the kernel and cokernel of $f$ are $v_{n}$-torsion.  
\item [(c)] For any comodule $M$, there is an exact sequence of comodules
\[
 0 \xrightarrow{} T_{n}M
   \xrightarrow{} M
   \xrightarrow{} L_{n}M
   \xrightarrow{} T'
   \xrightarrow{} 0
\]
where $T'$ is a $v_{n}$-torsion comodule. 
\item [(d)] A comodule $M$ is $L_{n}$-local if and only if 
\[
\Hom_{BP_{*}}^{*}(BP_{*}/I_{n+1},M) =
\Ext^{1,*}_{BP_{*}BP}(BP_{*}/I_{n+1},M)=0.
\]
\end{enumerate}
\end{proposition}
\begin{remark}
 We know from~\cite[Corollary 4]{landweber-comodules} that
 $\Hom_{BP_{*}}^{*}(BP_{*}/I_{n+1},M)=0$ iff
 $\Hom_{BP_{*}BP}^{*}(BP_{*}/I_{n+1},M)=0$.
\end{remark}
\begin{proof}
Parts~(a) and~(b) are immediate consequences of the fact that $L_{n}$ is
localization with respect to the $v_{n}$-torsion comodules.  Part~(d) is
proven in Corollary~4.3 of~\cite{hovey-strickland-comodules}.  Part~(c)
follows from the other parts; since $L_{n}$ is idempotent, the map
$\iota \mathcolon M\xrightarrow{} L_{n}M$ is an $L_{n}$-equivalence, so
its kernel and cokernel are $v_{n}$-torsion.  Part~(d) implies that
$L_{n}M$ has no $v_{n}$-torsion, so the kernel of $\iota$ is $T_{n}M$.
\end{proof}

In particular, of course, if $M$ is a $v_{n}$-torsion comodule, then
$L_{n}M=0$.  

We can use this proposition to identify some $L_{n}$-local comodules.  

\begin{proposition}\label{prop-local-regular}
Suppose $m<n$ and $M$ is a $v_{m-1}$-torsion comodule on which $(v_{m},
v_{m+1})$ is a regular sequence.  Then $M$ is $L_{n}$-local.
\end{proposition}

For this proposition, we need the following lemma.  We will prove the
converse of this lemma in Corollary~\ref{cor-local-Cech}.  

\begin{lemma}\label{lem-local-no-BP}
Suppose $M$ is a $BP_{*}BP$-comodule with no $v_{n}$-torsion such that 
\[
\Ext^{1,*}_{BP_{*}}(BP_{*}/I_{n+1},M) =0.
\]
Then $M$ is $L_{n}$-local.  
\end{lemma}

\begin{proof}
We must show that $\Ext^{1}_{BP_{*}BP}(BP_{*}/I_{n+1},M)=0$.  So suppose we
have a short exact sequence 
\[
0 \xrightarrow{} M \xrightarrow{} X \xrightarrow{p} BP_{*}/I_{n+1}
\xrightarrow{} 0
\]
of comodules.  We know that $X\cong M\bigoplus BP_{*}/I_{n+1}$ as
$BP_{*}$-modules, and $T_nM=0$, so $p$ induces an isomorphism
$T_nX\xrightarrow{}T_n(BP_*/I_{n+1})=BP_*/I_{n+1}$ (of comodules).
The inverse of this isomorphism splits the sequence.
\end{proof}

\begin{proof}[Proof of Proposition~\ref{prop-local-regular}] First note
that $M$ has no $v_{m}$-torsion, so $M$ has no $v_{n}$-torsion either.
In light of Lemma~\ref{lem-local-no-BP}, we must show that
$\Ext^{1,*}_{BP_{*}}(BP_{*}/I_{n+1},M)=0$.  For notational simplicity,
we will assume that $*=0$, but the proof works for any value of $*$.
Suppose we have a short exact sequence
\begin{equation}\label{eq-local}
0 \xrightarrow{} M \xrightarrow{f} X \xrightarrow{g} BP_{*}/I_{n+1}
\xrightarrow{} 0
\end{equation}
of $BP_{*}$-modules.  Choose an $x\in X$ such that $g(x)=1$.  Then
$g(v_{m}x)=g(v_{m+1}x)=0$, so there are elements $y$ and $z$ in $M$ such
that $f(y)=v_{m}x$ and $f(z)=v_{m+1}x$.  Then $f(v_{m+1}y)=f(v_{m}z)$,
so $v_{m+1}y=v_{m}z$.  Since $(v_{m}, v_{m+1})$ is a regular sequence on
$M$, we conclude that $y=v_{m}w$ for some $w$.  This means
$v_{m}v_{m+1}w=v_{m}z$, so, since $v_{m}$ is not a zero-divisor on $M$,
$z=v_{m+1}w$.  Now consider the element $x'=x-f(w)$.  We claim that this
element defines a splitting of the short exact sequence~\ref{eq-local}.
Certainly $g(x')=1$ and $v_{m}(x')=v_{m+1}(x')=0$.  Now
suppose $i\leq n$.  We claim that $v_{i}(x')=0$.  Certainly
$g(v_{i}(x'))=0$, so $v_{i}(x')=f(t)$ for some $t$.  But then
$f(v_{m}t)=0$, so $v_{m}t=0$.  This forces $t$ to be $0$, as required.
\end{proof}

Proposition~\ref{prop-local-regular} immediately gives us part of
Theorem~\ref{main-B}(5).  

\begin{corollary}\label{cor-local-easy-case}
Suppose $k<n$.  Then $BP_{*}/I_{k}$ is $L_{n}$-local.  
\end{corollary}

Another class of local comodules is given by the following proposition,
which proves part of Theorem~\ref{main-B}(4).  

\begin{proposition}\label{prop-local-periodic}
If $v_{m}$ acts invertibly on a $BP_{*}BP$-comodule $M$ for some $m\leq
n$, then $M$ is $L_{n}$-local.  
\end{proposition}

\begin{proof}
If $v_{m}$ acts invertibly on $M$ and $m<n$, then $(v_{m}, v_{m+1})$ is
a regular sequence on $M$, so Proposition~\ref{prop-local-regular}
implies $M$ is $L_{n}$-local.  If $v_{n}$ acts invertibly on $M$, then
certainly $M$ has no $v_{n}$-torsion.  By Lemma~\ref{lem-local-no-BP},
it suffices to show that $\Ext^{1,*}_{BP_{*}}(BP_{*}/I_{n+1},M)=0$.
Since $v_{n}$ acts by $0$ on $BP_{*}/I_{n+1}$, it also acts by $0$ on
this $\Ext^{1,*}$ group.  On the other hand, since $v_{n}$ acts
isomorphically on $M$, it acts isomorphically on this $\Ext^{1,*}$ group
as well.  Hence $\Ext^{1,*}(BP_{*}/I_{n+1},M)=0$ as required.  
\end{proof}

The following corollary proves part of Theorem~\ref{main-B}(6).

\begin{corollary}\label{cor-local-periodic}
Suppose $M$ is a $v_{n-1}$-torsion $BP_{*}BP$-comodule.  Then
$L_{n}M\cong v_{n}^{-1}M$.  
\end{corollary}

Note that this includes the case $n=0$, where we interpret $v_{0}=p$ and
$v_{-1}=0$.  

\begin{proof}
Note that $v_{n}^{-1}M$ is a $BP_{*}BP$-comodule
by~\cite[Proposition~2.9]{johnson-yosimura}.  The map
$M\xrightarrow{}v_{n}^{-1}M$ obviously has $v_{n}$-torsion kernel and
cokernel, so is an $L_{n}$-equivalence.
Proposition~\ref{prop-local-periodic} implies that $v_{n}^{-1}M$ is
$L_{n}$-local, so it must be $L_{n}M$.  
\end{proof}

Finally, we prove the part of Theorem~\ref{main-B}(8) dealing with
$L_{n}$ itself. 

\begin{proposition}\label{prop-local-filtered}
The functors $\Phi^{*}$ and $L_{n}$ commute with filtered colimits and
arbitrary direct sums.  
\end{proposition}

\begin{proof}
We show that $\Phi^{*}$ preserves filtered colimits.  It will then
follow that $\Phi^{*}$ preserves arbitrary direct sums, which are
filtered colimits of finite direct sums, and that
$L_{n}=\Phi^{*}\Phi_{*}$ preserves filtered colimits and arbitrary
direct sums, completing the proof.  

So suppose $X_{i}$ is a filtered diagram of $E(n)_{*}E(n)$-comodules.
There is certainly a natural map 
\[
\alpha \mathcolon \colim \Phi^{*}X_{i}\xrightarrow{}\Phi^{*}(\colim X_{i}).
\]
We need to recall that a $BP_{*}BP$-comodule
(resp. $E(n)_{*}E(n)$-comodule) $P$ is called \textbf{dualizable} if it
is finitely generated and projective over $BP_{*}$ (resp. $E(n)_{*}$),
and that the dualizable comodules generate the category of
$BP_{*}BP$-comodules (resp. $E(n)_{*}E(n)$-comodules).
See~\cite[Section~1.4]{hovey-comodule-homotopy}.  Thus, to show $\alpha
$ is an isomorphism, it suffices to check that it is an isomorphism upon
applying $BP_{*}BP \comod (P,-)$ for any dualizable $BP_{*}BP$-comodule
$P$, since the dualizable comodules generate.  The main point is that
$\Phi_{*}$ preserves dualizable comodules, and dualizable comodules are
finitely presented.  This implies
\begin{gather*}
BP_{*}BP \comod (P, \colim \Phi^{*}X_{i}) \cong \colim BP_{*}BP \comod (P,
\Phi^{*}X_{i}) \\
\cong \colim E(n)_{*}E(n)\comod (\Phi_{*}P, X_{i}) \cong
E(n)_{*}E(n)\comod (\Phi_{*}P, \colim X_{i}) \\ 
\cong BP_{*}BP \comod (P, \Phi^{*}(\colim X_{i})),
\end{gather*}
so $\Phi^{*}$ preserves filtered colimits. 
\end{proof}

\section{Injective $BP_{*}BP$-comodules}\label{sec-injective}

In order to construct the spectral sequence of Theorem~\ref{main-A} and
in order to compute the right derived functors of $L_{n}$, we need to
known something about injective objects in the category of
$BP_{*}BP$-comodules.  Very little seems to have been written about
these absolute injectives; relative injectives are easier to understand
and have been used much more often.  The object of this section is to
learn a little more; in particular, we prove that $L_{n}$, $\Phi_{*}$,
$\Phi^{*}$, and $T_{n}$ all preserve injectives.  

The most basic fact about injective $BP_{*}BP$-comodules is the
following well-known lemma.  Recall that a Hopf algebroid $(A, \Gamma)$
is said to be \textbf{flat} if $\Gamma$ is flat as a left (or,
equivalently, right) $A$-module.

\begin{lemma}\label{lem-enough-injectives}
Let $(A, \Gamma)$ be a flat Hopf algebroid. 
\begin{enumerate}
\item [(a)] If $I$ is an injective $A$-module, then the extended
$\Gamma$-comodule $\Gamma \otimes_{A}I$ is an injective
$\Gamma$-comodule.
\item [(b)] There are enough injective $\Gamma $-comodules. 
\item [(c)] A $\Gamma$-comodule is injective if and only if it is a
comodule retract of $\Gamma \otimes_{A}I$ for some injective $A$-module
$I$.
\end{enumerate}
\end{lemma}

\begin{proof}
Because the extended comodule functor is right adjoint to the forgetful
functor from $\Gamma$-comodules to $A$-modules, we have 
\[
\Hom_{\Gamma}(-,\Gamma \otimes_{A}M) \cong \Hom (-,M)
\]
from which part~(a) follows.  
  
Now, if $M$ is an arbitrary $\Gamma $-comodule, choose an injective
$A$-module $J$ so that there is an embedding $M\xrightarrow{j}J$.
The composite 
\[
M \xrightarrow{\psi } \Gamma \otimes_{A} M \xrightarrow{1\otimes
j} \Gamma \otimes_{A} J
\]
is a comodule embedding of $M$ into an injective $\Gamma $-comodule,
proving part~(b).  If $M$ is itself injective, this embedding must have a
retraction, proving part~(c).  
\end{proof}

This lemma is of little practical assistance, since injective
$BP_{*}$-modules are extremely complex.  They must not only be
$v_{n}$-divisible for all $n$, but also $x$-divisible for every
nonzero homogeneous element $x$ in $BP_{*}$.  This is the reason one
generally uses relatively injective $BP_{*}BP$-comodules, as they are
much simpler.  However, to compute right derived functors of $L_{n}$, we
must use absolute injectives.  

The first step is to understand the $v_{n}$-torsion in an injective
comodule.  

\begin{proposition}\label{prop-injective-torsion}
Suppose $M$ is a $v_{n}$-torsion $BP_{*}BP$-comodule and $N$ is an
essential extension of $M$ in the category of $BP_{*}BP$-comodules.
Then $N$ is $v_{n}$-torsion.  In particular, the injective hull of $M$
is $v_{n}$-torsion.
\end{proposition}

\begin{proof}
Suppose $N$ is not $v_{n}$-torsion.  Let $x$ be an element of $M$ that
is not $v_{n}$-torsion, and let $I=\sqrt{\Ann x}$.  Since $x$ is not
$v_{n}$-torsion, $v_{n}$ is not in $I$.  Theorem~1
of~\cite{landweber-comodules} guarantees that $I$ is an invariant ideal
of $BP_{*}$, so we must have $I=I_{k}$ for some $k\leq n$.  Theorem~2
of~\cite{landweber-comodules} tells us that there is some primitive $y$
in $N$ such that $\Ann (y)=I_{k}$.  Hence $BP_{*}/I_{k}$ is isomorphic
to a subcomodule of $N$.  This subcomodule has no $v_{n}$-torsion, so
cannot intersect $M$ nontrivially.  This contradicts our assumption that
$N$ is an essential extension of $M$.
\end{proof}

This proposition leads to the following useful theorem. 

\begin{theorem}\label{thm-injective-torsion}
Suppose $I$ is an injective $BP_{*}BP$-comodule, and let $T_{n}I$ denote
the $v_{n}$-torsion in $I$.  Then $T_{n}I$ and $I/T_{n}I$ are injective,
and $I\cong T_{n}I\oplus I/T_{n}I$.
\end{theorem}

\begin{proof}
The injective hull of $T_{n}I$ must be a subcomodule of $I$, since $I$
is injective, and it must be $v_{n}$-torsion by
Proposition~\ref{prop-injective-torsion}.  Hence it must be $T_{n}I$
itself.  
\end{proof}

\begin{corollary}\label{cor-injective-Ln}
Suppose $I$ is an injective $BP_{*}BP$-comodule.  Then
$L_{n}I=I/T_{n}I$.  In particular, $L_{n}$ preserves injectives.  
\end{corollary}

\begin{proof}
Certainly the map $I\xrightarrow{}I/T_{n}I$ is an $L_{n}$-equivalence.
But $I/T_{n}I$ is is an injective by
Theorem~\ref{thm-injective-torsion}, and has no $v_{n}$-torsion, so is
$L_{n}$-local.  
\end{proof}

\begin{corollary}\label{cor-injective-phi}
The functor $\Phi^{*}$ preserves and reflects injectivity, and the
functor $\Phi_{*}$ preserves injectives. 
\end{corollary}

\begin{proof}
The functor $\Phi^{*}$ is right adjoint to the exact functor $\Phi_{*}$,
so preserves injectives.  Conversely, suppose $\Phi^{*}I$ is injective,
$j\mathcolon M\xrightarrow{j}N$ is an inclusion of
$E(n)_{*}E(n)$-comodules, and $f\mathcolon M\xrightarrow{}I$ is a map.
Applying $\Phi^{*}$, we find a map $h\mathcolon
\Phi^{*}N\xrightarrow{}\Phi^{*}I$ such that $h\circ
\Phi^{*}j=\Phi^{*}f$.  Since $\Phi^{*}$ is fully faithful, we conclude
that $h=\Phi^{*}g$ for some extension $g$ of $f$.  Hence $I$ is
injective.  

Now $L_{n}=\Phi^{*}\Phi_{*}$ preserves injectives by
Corollary~\ref{cor-injective-Ln}.  Since $\Phi^{*}$ reflects injectives,
we conclude that $\Phi_{*}$ must preserve injectives.  
\end{proof}

Theorem~\ref{thm-injective-torsion} divides the study of injective
$BP_{*}BP$-comodules into those with no $v_{n}$-torsion and those which
are all $v_{n}$-torsion.  About all we know about injective comodules
which are all $v_{n}$-torsion is the following proposition.  

\begin{proposition}\label{prop-injective-divisible}
Suppose $I$ is an injective $BP_{*}BP$-comodule that is all
$v_{n}$-torsion.  Then $I$ is $v_{n+1}$-divisible. 
\end{proposition}

\begin{proof}
Suppose $x\in I$.  Because every $BP_{*}BP$-comodule is a filtered
colimit of finitely presented $BP_{*}BP$-comodules, there is a map
$P\xrightarrow{}I$ from a comodule $P$ that is a free finitely generated
$BP_{*}$-module, whose image contains $x$.  Since $I$ is
$v_{n}$-torsion, and therefore $v_{i}$-torsion for all $i\leq n$, this
map factors through
\[
Q=P/JP \xrightarrow{g} I
\]
for some invariant ideal $J=(p^{i_{0}}, v_{1}^{i_{1}},\dotsc
,v_{n-1}^{i_{n-1}})$.  There is some $k$ such that $v_{n+1}^{k}$ is
invariant modulo $J$.  Thus multiplication by $v_{n+1}^{k}$ defines a
monomorphism of comodules $Q\xrightarrow{}\Sigma^{-t}Q$.  Because $I$
is injective, $g$ must extend to a map $\Sigma^{-t}Q\xrightarrow{}I$,
showing that $x$ is divisible by $v_{n+1}$.
\end{proof}

We now turn our attention to injectives that have no $v_{n}$-torsion. 

\begin{theorem}\label{thm-injective-no-torsion}
Suppose $M$ is a $BP_{*}BP$-comodule with no $v_{n}$-torsion.  Then
there is an embedding of $M$ into an injective $BP_{*}BP$-comodule with
no $v_{n}$-torsion, and this embedding can be chosen to be functorial on
the category of $BP_{*}BP$-comodules with no $v_{n}$-torsion.  
\end{theorem}

\begin{proof}
Since the category of $E(n)_{*}E(n)$-comodules is a Grothendieck
category (see Section~1.4 of~\cite{hovey-comodule-homotopy}),
there is a functorial embedding of any $E(n)_{*}E(n)$-comodule into an
injective $E(n)_{*}E(n)$-comodule.  (Apply Quillen's small object
argument to the set of subobjects of a generator).  In particular, for
$M$ a $BP_{*}BP$-comodule, we get a functorial embedding
$\Phi_{*}M\xrightarrow{}I$.  Applying $\Phi^{*}$ gives us a functorial
embedding $L_{n}M\xrightarrow{}\Phi^{*}I$, and $\Phi^{*}I$ is an
injective comodule (by Corollary~\ref{cor-injective-phi}), and has no
$v_{n}$-torsion (since it is $L_{n}$-local).  Since $M$ has no
$v_{n}$-torsion, $M$ embeds in $L_{n}M$.  
\end{proof}

We expect that injective $BP_{*}BP$-comodules are not closed under
filtered colimits, though we do not have a counterexample.  Those with
no $v_{n}$-torsion, on the other hand, are better behaved.  

\begin{proposition}\label{prop-injective-filtered}
Injective $BP_{*}BP$-comodules with no $v_{n}$-torsion are closed under
filtered colimits.  
\end{proposition}

This proposition depends on the following lemma.

\begin{lemma}\label{lem-injective-Adams}
Injective $E(n)_{*}E(n)$-comodules are closed under filtered colimits. 
\end{lemma}

\begin{proof}
Recall that the category of $E(n)_{*}E(n)$-comodules is a Grothendieck
category; a set of generators is given by the comodules which are
finitely generated and projective over
$E(n)_{*}$~\cite[Section~1.4]{hovey-comodule-homotopy}.  There is a
version of Baer's criterion for injectivity that works for any
Grothendieck category~\cite{stenstrom}.  Let $\{G_{j} \}$ be a set of
generators for the Grothendieck category in question; then an object $I$
is injective if and only if $\Hom (G_{j},I)\xrightarrow{}\Hom (N_{j},I)$
is surjective for all $j$ and all subobjects $N_{j}$ of $G_{j}$.  In
particular, if $\Hom (G_{j},-)$ and $\Hom (N_{j},-)$ commute with
filtered colimits (that is, if $G_{j}$ and $N_{j}$ are finitely
presented), then injectives are closed under filtered colimits.  In our
case, the generators $G_{j}$ are finitely generated and projective over
$E(n)_{*}$.  Since $E(n)_{*}$ is Noetherian, the objects $N_{j}$ are
also finitely presented over $E(n)_{*}$.  This means the objects $N_{j}$
and $G_{j}$ are also finitely presented as $E(n)_{*}E(n)$-comodules by
Proposition~1.3.3 of~\cite{hovey-comodule-homotopy}, completing the
proof.
\end{proof}

\begin{proof}[Proof of Proposition~\ref{prop-injective-filtered}]
Suppose $F\mathcolon \cat{J}\xrightarrow{}BP_{*}BP\comod$ is a functor
from a filtered category $\cat{J}$ to injective comodules with no
$v_{n}$-torsion.  Then $F(j)$ is $L_{n}$-local for all $j\in J$, so we
have 
\[
\colim F \cong \colim \Phi^{*}\Phi_{*}F \cong \Phi^{*}(\colim
\Phi_{*}F). 
\]
Now $\Phi_{*}F(j)$ is an injective $E(n)_{*}E(n)$-comodule for all $j\in
\cat{J}$ by Corollary~\ref{cor-injective-phi}, so
Lemma~\ref{lem-injective-Adams} tells us that $\colim \Phi_{*}F$
is injective.  Since $\Phi^{*}$ preserves injectives, we conclude that
$\colim F$ is injective. 
\end{proof}

We can now give a partial structure theorem for injective
$BP_{*}BP$-comodules.

\begin{proposition}\label{prop-injective-decomp}
Suppose $I$ is an injective $BP_{*}BP$-comodule, and $n\geq 0$.  Then
\[
I \cong I_{0} \oplus I_{1} \oplus \dotsb I_{n} \oplus T_nI
\]
where 
\begin{itemize}
 \item[(a)] Each $I_j$ is an injective $BP_{*}BP$-comodule with
  $v_j^{-1}I_j=I_j$ \ulp and thus $I_j$ is $v_{j-1}$-torsion\urp .
 \item[(b)] $T_nI$ is injective and $v_n$-torsion.
\end{itemize}
In particular, if $I$ is indecomoposable then either $I=v_{j}^{-1}I$
for some $j$ or $I$ is $v_{j}$-torsion for all $j$.
\end{proposition}

\begin{proof}
Put $I_j=T_{j-1}I/T_jI$ (where $T_{-1}I=I$).  As the comodules $T_jI$
are injective (by Theorem~\ref{thm-injective-torsion}), the filtration 
\[ T_nI\leq T_{n-1}I\leq\ldots\leq T_0I\leq I \]
must split, giving $I=T_nI\oplus\bigoplus_{j=0}^{n}I_j$.  The comodule
$I_j$ is a summand of $I$ and thus is injective.  By construction it
is $v_{j-1}$-torsion, and thus $v_j$-divisible by
Proposition~\ref{prop-injective-divisible}.  The definition also
implies that there is no $v_j$-torsion, so $I_j=v_j^{-1}I_j$.
\end{proof}

It would be nice to have some explicit knowledge of injective
$BP_{*}BP$-comodules $I$ with $v_{n}^{-1}I=I$.  When $n=0$, at least,
this is easy. 

\begin{proposition}\label{prop-injective-rational}
Suppose $M$ is a $BP_{*}BP$-comodule with no $p$-torsion.  Then $M$ is
injective if and only if $M$ is a rational vector space.  
\end{proposition}

\begin{proof}
Proposition~\ref{prop-injective-divisible} shows that if $M$ is
injective, then it must be a rational vector space.  Conversely, if $M$
is rational, then $M=L_{0}M=\Phi^{*}\Phi_{*}M$.  The category of
$E(0)_{*}E(0)$-comodules is the category of rational vector spaces, so
$\Phi_{*}M$ is injective.  Since $\Phi^{*}$ preserves injectives, we
conclude that $M$ is injective.  
\end{proof}

The analogue of this proposition is definitely false when $n>0$.  Still,
this gives a rationale for why the chromatic resolution is useful.
Indeed, suppose we want to find an injective resolution of $BP_{*}$ as a
$BP_{*}BP$-comodule.  Proposition~\ref{prop-injective-rational} implies
that $M^{0}=p^{-1}BP_{*}$ is the injective hull of $BP_{*}$ as a
$BP_{*}BP$-comodule.  The cokernel $N^{1}$ is usually written
$BP_{*}/(p^{\infty})$.  The injective hull of $N^{1}$ must be a
$p$-torsion essential extension of $N^{1}$ on which $v_{1}$ acts
invertibly.  The simplest way to do this is to form
$M^{1}=v_{1}^{-1}N^{1}$, which is the next term in the chromatic
resolution.  Sadly, $N^{1}$ is not actually injective, but it seems to
be the closest one can get to the injective hull of $M^{1}$ in a fairly
simple way.  Iterating this idea leads to the chromatic resolution.  

\section{The derived functors of $L_{n}$}\label{sec-derived}

Now that we have some knowledge of injective $BP_{*}BP$-comodules, we
can begin to compute derived functors.  The goal of this section is to
complete the proof of Theorem~\ref{main-B} except for part~(1), which we
deal with in the next section.  

Recall that $L^{i}_{n}$ denotes the $i$th right derived functor of
$L_{n}$.  We also let $T^{i}_{n}$ denote the $i$th right derived functor
of $T_{n}$, where $T_{n}(M)$ is the subcomodule of $v_{n}$-torsion
elements in $M$.

The first thing to point out is that $L^{i}_{n}$ and $T^{i}_{n}$ are
closely related. 

\begin{theorem}\label{thm-derived-different}
For $M$ a $BP_{*}BP$-comodule, we have a natural short exact sequence 
\[
0 \xrightarrow{} T_{n}M \xrightarrow{} M \xrightarrow{} L_{n}M
\xrightarrow{} T^{1}_{n}M \xrightarrow{} 0.  
\]
and natural isomorphisms 
\[
L^{i}_{n}M \cong T^{i+1}_{n}M
\]
for $i>0$.  
\end{theorem}

\begin{proof}
Let $I_{*}$ be an injective resolution of $M$.  Then
$T^{i}_{n}M\cong H_{-i}(T_{n}I_{*})$, and $L^{i}_{n}M\cong
H_{-i}(L_{n}I_{*})$.  But $L_{n}I_{*}\cong I_{*}/T_{n}I_{*}$ by
Corollary~\ref{cor-injective-Ln}.  Hence we have a natural short exact
sequence of complexes 
\[
0 \xrightarrow{} T_{n}I_{*} \xrightarrow{} I_{*} \xrightarrow{}
L_{n}I_{*} \xrightarrow{} 0.  
\]
The long exact sequence in homology gives the desired result.  
\end{proof}

We also point out that computing $L^{i}_{n}$ is equivalent to computing
the right derived functors of $\Phi^{*}$.  

\begin{proposition}\label{prop-derived-phi}
Let $M$ be a $BP_{*}BP$-comodule, and let $R^{i}\Phi^{*}$ denote the
$i$th right derived functor of $\Phi^{*}$.  Then we have a natural
isomorphism 
\[
(R^{i}\Phi^{*})(\Phi_{*}M) \cong L^{i}_{n}M.
\]
\end{proposition}

Note that, since $\Phi_{*}\Phi^{*}N\cong N$, we can also write this
isomorphism as 
\[
(R^{i}\Phi^{*})(N) \cong L^{i}_{n}(\Phi^{*}N).  
\] 

\begin{proof}
Let $I_{*}$ be an injective resolution of $M$.  Then $\Phi_{*}I_{*}$ is
an injective resolution of $\Phi_{*}M$, since $\Phi_{*}$ is exact and
preserves injectives by Corollary~\ref{cor-injective-phi}.  Hence 
\[
(R^{i}\Phi^{*})(\Phi_{*}M) \cong H_{-i}(\Phi^{*}\Phi_{*}I_{*})\cong
H_{-i}(L_{n}I_{*}) \cong L^{i}_{n}M.  
\]
\end{proof}

We now begin the computation of $L^{i}_{n}$.  

\begin{proposition}\label{prop-torsion-no-derived}
If $T$ is a $v_{n}$-torsion $BP_{*}BP$-comodule, then $L^{i}_{n}T=0$ for
all $i\geq 0$.  Furthermore, for an arbitrary comodule $M$, the map
$M\xrightarrow{}L_{n}M$ induces an isomorphism
$L^{i}_{n}M\xrightarrow{}L^{i}_{n}L_{n}M$.
\end{proposition}

This proposition proves part~(7) of Theorem~\ref{main-B}.  

\begin{proof}
Using Proposition~\ref{prop-injective-torsion}, one can easily construct
an injective resolution $I_{*}$ of $T$ that is all $v_{n}$-torsion.
Hence $L_{n}I_{*}=0$, so $L^{i}_{n}T=0$ for all $i$.  For the second
statement, recall that we have short exact sequences 
\[
0 \xrightarrow{} T \xrightarrow{} M \xrightarrow{} M/T \xrightarrow{} 0
\]
and 
\[
0 \xrightarrow{} M/T \xrightarrow{} L_{n}M \xrightarrow{} T'
\xrightarrow{} 0
\]
where $T$ and $T'$ are $v_{n}$-torsion.  Applying $L_{n}$ gives the
desired result.  
\end{proof}

The following proposition is part~(3) of Theorem~\ref{main-B}.  

\begin{proposition}\label{prop-derived-are-torsion}
Suppose $M$ is a $BP_{*}BP$-comodule.  Then $L^{i}_{n}M$ is
$v_{n}$-torsion for $i>0$.  
\end{proposition}

\begin{proof}
Let $I_{*}$ be an injective resolution of $M$.  Then
$T^{i}_{n}(M)=H_{-i}T_{n}I_{*}$ is obviously $v_{n}$-torsion.  The
result follows from Theorem~\ref{thm-derived-different}.  
\end{proof}

We now show that the chromatic resolution is as good as an injective
resolution for computing $L^{i}_{n}$.  This completes the proof of
part~(4) of Theorem~\ref{main-B}. 

\begin{theorem}\label{thm-chromatic-derived}
 Suppose $M$ is a $BP_{*}BP$-comodule on which $v_{k}$ acts
 isomorphically for some $k$.  Then $L^{i}_{n}M=0$ for all $i>0$ and
 all $n$.  Moreover, we have $L_nM=0$ if $n<k$ and $L_nM=M$ if 
 $n\geq k$.
\end{theorem}

\begin{proof}
We claim that we can choose an injective resolution $I_{*}$ of $M$ for
which $I_{*}=v_{k}^{-1}I_{*}$.  To see this, it suffices by induction to
show that if $N$ is a $BP_{*}BP$-comodule for which $N=v_{k}^{-1}N$,
there there is a short exact sequence 
\[
0 \xrightarrow{} N \xrightarrow{} I \xrightarrow{} N' \xrightarrow{} 0
\]
of comodules for which $I$ is injective, $v_{k}^{-1}I=I$, and
$v_{k}^{-1}N'=N'$.  Since $N=v_{k}^{-1}N$, $N$ is all $v_{k-1}$-torsion
by Proposition~2.9 of~\cite{johnson-yosimura}, and of course $N$ has no
$v_{k}$-torsion.  Proposition~\ref{prop-injective-torsion} and
Theorem~\ref{thm-injective-no-torsion} together imply that the injective
hull $I$ of $N$ is $v_{k-1}$-torsion and has no $v_{k}$-torsion.
Proposition~\ref{prop-injective-divisible} then implies that
$I=v_{k}^{-1}I$.  It follows easily that multiplication by $v_{k}$ is
surjective on $N'$, but we claim it is injective as well.  Indeed,
suppose $x\in N'$ has $v_{k}x=0$.  Choose a $y$ in $I$ whose image in
$N'$ is $x$, so that $v_{k}y$ is in $N'$.  Since $N'=v_{k}^{-1}N'$,
there is a $z$ in $N'$ such that $v_{k}z=v_{k}y$.  It follows that
$z=y$, and so $x=0$.  

We now have our desired injective resolution $I_{*}$ of $M$ for which
$I_{*}=v_{k}^{-1}I_{*}$.  The argument now breaks into two cases.  If
$n\geq k$, we apply the $v_{n}$-torsion functor $T_{n}$.  Since there is
no $v_{k}$-torsion in $I_{*}$, there is also no $v_{n}$-torsion by
Lemma~2.3 of~\cite{johnson-yosimura}.  Thus $T_{n}I_{*}=0$, and so 
\[
L^{i}_{n}M\cong T^{i+1}_{n}M \cong H_{-i-1}T_{n}I_{*}=0
\]
for all $i>0$.  For $i=0$, use Proposition~\ref{prop-local-periodic}.

Now suppose $n<k$.  Since $v_{k}^{-1}I_{*}=I_{*}$, $I_{*}$ is all
$v_{k-1}$-torsion, so also all $v_{n}$-torsion.  Hence $L_{n}I_{*}=0$,
so 
\[
L^{i}_{n}M =H_{-i}L_{n}I_{*}=0
\]
for $i\geq 0$. 
\end{proof}

Theorem~\ref{thm-chromatic-derived} allows us to compute $L^{i}_{n}M$
some important $BP_{*}BP$-comodules $M$.   This completes the proof of
parts~(5) and~(6) of Theorem~\ref{main-B}.  

\begin{corollary}\label{cor-derived-standard}
\begin{enumerate}
\item [(a)] Suppose $n<k$.  Then $L^{i}_{n}(BP_{*}/I_{k})=0$ for all
$i$.  
\item [(b)] $L^{i}_{k}(BP_{*}/I_{k})=0$ for $i>0$, whereas
$L^{0}_{k}(BP_{*}/I_{k})=v_{k}^{-1}BP_{*}/I_{k}$.  
\item [(c)] Suppose $n>k$.  Then $L^{i}_{n}(BP_{*}/I_{k})=0$ unless
$i=0$ or $n-k$.  We have 
\[
L^{0}_{k}(BP_{*}/I_{k})=BP_{*}/I_{k}
\]
and
\[
L^{n-k}_{k}(BP_{*}/I_{k})=BP_{*}/(p,v_{1},\dotsc ,v_{k-1},
v_{k}^{\infty},\dotsc ,v_{n}^{\infty}).
\]
\end{enumerate}
\end{corollary}

\begin{proof}
Let $M=BP_{*}/I_{k}$, and consider the chromatic resolution
$M\xrightarrow{}J_{*}$ of $M$, where
\[
J_{t}=v_{t+k}^{-1}BP_{*}/(p,v_{1},\dotsc ,v_{k-1},v_{k}^{\infty},\dotsc
,v_{t+k-1}^{\infty}).  
\]
By Theorem~\ref{thm-chromatic-derived}, we have
$L^{i}_{n}J_{t}=0$ for all $i>0$.  Hence $L^{i}_{n}M\cong
H_{-i}L_{n}J_{*}$.  Now, each of the comodules $J_{k}$ is
$v_{k-1}$-torsion, so $L_{n}J_{*}=0$ if $n<k$.  This completes the proof
of part~(a).  

If $n=k$, then $L_{n}J_{t}=0$ for $t>0$, from which part~(b) follows
easily.  If $n>k$, on the other hand, $L_{n}J_{t}=J_{t}$ for $t<n-k+1$,
from which part~(c) follows.  
\end{proof}

We also discover that $L_{n}$ has only finitely many right derived
functors.  

\begin{theorem}\label{thm-torsion-some-derived}
Suppose $M$ is a $v_{k}$-torsion comodule for some $-1\leq k\leq n$.
Then $L^{i}_{n}M=0$ for $i\geq n-k$.  In particular, $L^{i}_{n}N=0$ for
$i>n$ for any comodule $N$.  
\end{theorem}

For the purposes of this theorem, we take $v_{-1}=0$, so that every
comodule is $v_{-1}$-torsion.  This theorem proves part~(2) of
Theorem~\ref{main-B}.  

\begin{proof}
We proceed by downwards induction on $k$.  The base case $k=n$ is
Proposition~\ref{prop-torsion-no-derived}.  So suppose we know the theorem
for $k$, and $M$ is a $v_{k-1}$-torsion comodule.  Let $T_{k}M$ denote the
$v_{k}$-torsion in $M$.  We have a short exact sequence 
\[
0 \xrightarrow{} T_{k}M \xrightarrow{} M \xrightarrow{} N \xrightarrow{} 0
\]
where $N$ has no $v_{k}$-torsion.  By our induction hypothesis,
$L^{i}_{n}(T_{k}M)=0$ for $i\geq n-k$.  It therefore suffices to show
that $L^{i}_{n}(N)=0$ for $i>n-k$.

Now, since $N$ is $v_{k-1}$-torsion but has no $v_{k}$-torsion, we have
a short exact sequence 
\[
0 \xrightarrow{} N \xrightarrow{} v_{k}^{-1}N \xrightarrow{} T
\xrightarrow{} 0, 
\]
where $T$ is $v_{k}$-torsion.  Our induction hypothesis guarantees that
$L^{i}_{n}T=0$ for $i\geq n-k$, and Theorem~\ref{thm-chromatic-derived}
guarantees that $L^{i}_{n}T\cong L^{i+1}_{n}N$ for $i>0$.  Hence
$L^{i}_{n}N=0$ for $i>n-k$, as required.  
\end{proof}

Corollary~\ref{cor-derived-standard} together with the Landweber
filtration theorem gives a method for computing $L^{i}_{n}M$ for
finitely presented $BP_{*}BP$-comodules $M$.  To compute for more
general comodules $M$, we use the following theorem, which is part~(8)
of Theorem~\ref{main-B}.  

\begin{theorem}\label{thm-derived-filtered-colimits}
The functors $L^{k}_{n}$ preserve filtered colimits of
$BP_{*}BP$-comodules.  
\end{theorem}

Since $L_{n}$ itself preserves filtered colimits, this theorem would be
easy if filtered colimits of injective comodules were injective, but we
believe that this is false in general.  However, to compute $L^{i}_{n}$
the only injectives that matter are injectives with no $v_{n}$-torsion,
and these we know are closed under filtered colimits by
Proposition~\ref{prop-injective-filtered}.

\begin{proof}
We use induction on $k$.  When $k=0$ we have seen this already in
Proposition~\ref{prop-local-filtered}.  
Now suppose $L^{k}_{n}$ preserves filtered colimits for some $k\geq 0$,
and let $\{ M_{t}\}$ be a filtered diagram of comodules.  Then
$\{L_{n}M_{t}\}$ is a filtered diagram of comodules with no
$v_{n}$-torsion, so we can use Theorem~\ref{thm-injective-no-torsion} to
find a filtered diagram of injectives $\{I_{t} \}$ with no
$v_{n}$-torsion and a short exact sequence of filtered diagrams
\[
\{0 \} \xrightarrow{} \{L_{n}M_{t} \} \xrightarrow{}\{I_{t} \}
\xrightarrow{}\{N_{t} \} \xrightarrow{} \{0 \}.
\]
This gives us a short exact sequence 
\[
0 \xrightarrow{} \colim L_{n}M_{t} \xrightarrow{} \colim I_{t}
\xrightarrow{} \colim N_{t} \xrightarrow{} 0,
\]
and $\colim I_{t}$ is injective by
Proposition~\ref{prop-injective-filtered}.  

We must now separate the case $k=0$ from the case $k>0$.  If $k=0$, we
get exact sequences
\[
0 \xrightarrow{} \colim L_{n}M_{t} \xrightarrow{} \colim I_{t}
\xrightarrow{} \colim L_{n}N_{t} \xrightarrow{} \colim L^{1}_{n}M_{t}
\xrightarrow{} 0, 
\]
and 
\[
0 \xrightarrow{} \colim L_{n}M_{t} \xrightarrow{} \colim I_{t}
\xrightarrow{} L_{n}(\colim N_{t}) \xrightarrow{} L^{1}_{n}(\colim
L_{n}M_{t}) \xrightarrow{} 0.  
\]
There is a map from the first of these sequences to the second, which is
an isomorphism on every nonzero term except the last one, so we get an
isomorphism 
\[
\colim L^{1}_{n}M_{t} \cong L^{1}_{n}(\colim L_{n}M_{t}).  
\]
On the other hand, using Proposition~\ref{prop-torsion-no-derived}, and
the fact that $L_{n}$ commutes with filtered colimits, we get 
\[
L^{1}_{n}(\colim L_{n}M_{t}) \cong L^{1}_{n}L_{n}(\colim M_{t}) \cong
L^{1}_{n}(\colim M_{t}), 
\]
as required.  

If $k>0$, the situation is easier.  Indeed, using
Proposition~\ref{prop-torsion-no-derived} and the fact that $L_{n}$
commutes with filtered colimits, we have
\begin{gather*}
\colim L^{k+1}_{n}M_{t} \cong \colim L^{k+1}_{n}(L_{n}M_{t}) \cong
\colim L^{k}_{n}N_{t} \\
\cong L^{k}_{n}(\colim N_{t}) \cong
L^{k+1}_{n}(\colim L_{n}M_{t}) \cong L^{k+1}_{n}L_{n}(\colim M_{t})
\cong L^{k+1}_{n}(\colim M_{t}), 
\end{gather*}
completing the proof. 
\end{proof}

\section{Comparison with \v{C}ech cohomology}\label{sec-Cech}

The object of this section is to prove part~(1) of
Theorem~\ref{main-B}, showing that, for a comodule $M$, $L_{n}^{i}(M)$
is the same as the $i$th \textCech cohomology group
$\CechH_{I_{n+1}}^{i}M$ of $M$ with respect to $I_{n+1}$.  We also
show that \textCech cohomology $\CechH_{I_{n+1}}^{*}(-)$ is the
derived functors of of localization in the category of
$BP_{*}$-modules with respect to the hereditary torsion theory of
$I_{n+1}$-torsion modules.

We first remind the reader of the definition of \v{C}ech cohomology
from~\cite{greenlees-may-completions}.  Given an element $\alpha$ in a
commutative ring $R$, which we will always take to be $BP_{*}$, we form
the cochain complex $K^{\bt}(\alpha)$ which is $R$ in degree $0$ and
$R[1/\alpha]$ in degree $1$, with the differential being the obvious
map.  Given an ideal $I=(\alpha_{0},\dotsc ,\alpha_{n})$, we define
$K^{\bt}(I)$ to be the cochain complex
\[
K^{\bt}(I) = K^{\bt}(\alpha_{0}) \otimes_{R} K^{\bt}(\alpha_{1})
\otimes_{R} \dotsb \otimes_{R} K^{\bt}(\alpha_{n}).  
\]
This stable Koszul complex of course depends on the choice of generators
$\alpha_{i}$, but its quasi-isomorphism class does
not~\cite[Corollary~1.2]{greenlees-may-completions}.  There is an
obvious surjection $K^{\bt}(\alpha_{i})\xrightarrow{}R$ of complexes,
where $R$ is the complex consisting of $R$ concentrated in degree $0$.
Tensoring these together gives us a map 
\[
\epsilon \mathcolon K^{\bt}(I)\xrightarrow{}R. 
\]
We define the flat \textCech complex $\Cech ^{\bt}(I)$ by 
\[
\Sigma^{-1}\Cech^{\bt}(I) = \ker \epsilon.  
\]
Thus 
\[
\Cech^{k}(I) = \bigoplus_{|S|=k+1} R[1/\alpha_{S}]
\]
for $0\leq k\leq n$, where $S$ runs through the $k+1$-element subsets of
$(0,1,\dotsc ,n)$ and $\alpha_{S}=\prod_{i\in S}\alpha_{i}$.  

\begin{definition}\label{defn-local-cohomology}
The \textbf{local cohomology} $H_{I}^{*}(M)$ of an $R$-module $M$ with
respect to a finitely generated ideal $I=(\alpha_{0},\dotsc
,\alpha_{n})$ is
\[
H_{I}^{*}(M) = H^{*}(K^{\bt}(I)\otimes_{R}M).  
\]
The \textbf{\textCech cohomology} $\CechH_{I}^{*}(M)$ of $M$ with
respect to $I$ is 
\[
\CechH_{I}^{*}(M) = H^{*}(\Cech^{\bt}(I)\otimes_{R}M).  
\]
\end{definition}

Some of the basic properties of local and Cech cohomology are
summarized in the following proposition. 

\begin{proposition}\label{prop-Cech}
Suppose $I=(\alpha_{0},\dotsc ,\alpha_{n})$ is a finitely generated
ideal in a commutative ring $R$, and $M$ is an $R$-module.  
\begin{enumerate}
\item [(a)] We have a natural exact sequence 
\[
0 \xrightarrow{} H_{I}^{0}(M) \xrightarrow{} M \xrightarrow{}
\CechH_{I}^{0}(M) \xrightarrow{} H_{I}^{1}(M) \xrightarrow{} 0,
\]
and natural isomorphisms $\CechH_{I}^{k}(M)\cong H_{I}^{k+1}(M)$ for $k>0$.  
\item [(b)] $\CechH_{I}^{k}(M)=0$ unless $0\leq k\leq n$.  
\item [(c)] $H_{I}^{k}(M)$ is $I$-torsion for all $k$, and $\CechH_{I}^{k}(M)$ is $I$-torsion for all $k>0$.  On the other hand, $\CechH_{I}^{0}(M)$ has no $I$-torsion.  
\item [(d)] $H_{I}^{0}(M)$ is the submodule of $I$-torsion elements in
$M$.  
\item [(e)] $\CechH_{I}^{k}(M)=0$ for all $k$ if and only if $M$ is
$I$-torsion, and this is true if and only if $\CechH_{I}^{0}(M)=0$. 
\item [(f)] A short exact sequence 
\[
0 \xrightarrow{} M' \xrightarrow{} M \xrightarrow{} M'' \xrightarrow{} 0
\]
of $R$-modules gives rise to natural long exact sequences 
\begin{multline*}
0 \xrightarrow{} H_{I}^{0}(M') \xrightarrow{} H_{I}^{0}(M)
\xrightarrow{} H_{I}^{0}(M'') \\
\xrightarrow{} H_{I}^{1}(M') \xrightarrow{} \dotsb \xrightarrow{}
H_{I}^{n+1}(M) \xrightarrow{} H_{I}^{n+1}(M'') \xrightarrow{} 0, 
\end{multline*}
and 
\begin{multline*}
0 \xrightarrow{} \CechH_{I}^{0}(M') \xrightarrow{} \CechH_{I}^{0}(M)
\xrightarrow{} \CechH_{I}^{0}(M'') \\
\xrightarrow{} \CechH_{I}^{1}(M') \xrightarrow{} \dotsb \xrightarrow{}
\CechH_{I}^{n}(M) \xrightarrow{} \CechH_{I}^{n}(M'') \xrightarrow{} 0.
\end{multline*}
\item [(g)] Both $H_{I}^{0}$ and $\CechH_{I}^{0}$ are left exact
idempotent functors.
\end{enumerate}
\end{proposition}

Note that part~(f) certainly suggests that $H_{I}^{k}$ is the $k$th
right derived functor of $H_{I}^{0}$, but it does not prove it, since we
also need to know that $H_{I}^{k}$ sends injective modules to $0$ for
all $k>0$.  We will have to deal with this issue later.  

\begin{proof}
Most of this proposition follows from~\cite{greenlees-may-completions}.
Part~(a) appears in Section~1 of that paper, and part~(b) is obvious from
the definition of $\Cech^{\bt}(I)$.  The first sentence of part~(c) is
also in Section~1 of~\cite{greenlees-may-completions}.  For the second
part of part~(c), simply note that $\CechH_{I}^{0}(M)$ is a submodule
of $\bigoplus_{i}M[1/\alpha_{i}]$.  As $M[1/\alpha_{i}]$ has no
$\alpha_{i}$-torsion, it follows that $\bigoplus_{i}M[1/\alpha_{i}]$ has
no $I$-torsion.  Part~(d) is clear from the fact that $H_{I}^{0}(M)$ is
the kernel of the map 
\[
M \xrightarrow{} \bigoplus_{i} M[1/\alpha_{i}].  
\]

For part~(e), first suppose that $M$ is $I$-torsion.  Then the complex
$\Cech^{\bt}(I)\otimes_{R}M$ is the zero complex, and so of course
$\CechH_{I}^{k}(M)=0$ for all $k$.  Conversely, suppose $\CechH_{I}^{0}(M)=0$.  Then parts~(a) and~(d) show that $M$ is $I$-torsion.
For part~(f), simply note that the complexes $K^{\bt}(I)$ and
$\Cech^{\bt}(I)$ are complexes of flat modules.  Hence, a short exact
sequence of modules gives rise to a short exact sequence of complexes on
applying either $K^{\bt}(I)\otimes_{R}(-)$ or
$\Cech^{\bt}(I)\otimes_{R}(-)$.  The resulting long exact sequence in
cohomology gives us part~(f).  For part~(g), note that part~(f)
immediately implies that both $H_{I}^{0}$ and $\CechH_{I}^{0}$ are left
exact.  Part~(d) shows that $H_{I}^{0}$ is idempotent.  To see that
$\CechH_{I}^{0}$ is idempotent, apply part~(f) to the short exact
sequences
\[
0 \xrightarrow{} H_{I}^{0}(M) \xrightarrow{} M \xrightarrow{}
M/H_{I}^{0}(M) \xrightarrow{} 0
\]
and 
\[
0\xrightarrow{} M/H_{I}^{0}(M) \xrightarrow{}\CechH_{I}^{0}(M)
\xrightarrow{} H_{I}^{1}(M) \xrightarrow{} 0.  
\]
Since both $H_{I}^{0}(M)$ and $H_{I}^{1}(M)$ are $I$-torsion by
part~(c), part~(e) tells us that we get isomorphisms 
\[
\CechH_{I}^{0}(M) \cong \CechH_{I}^{0}(M/H_{I}^{0}(M)) \cong \CechH_{I}^{0}(\CechH_{I}^{0}M),
\]
completing the proof.  
\end{proof}

\begin{corollary}\label{cor-Cech-local}
Suppose $I=(\alpha_{0},\alpha_{1},\dotsc ,\alpha_{n})$ is a finitely
generated ideal in a commutative ring $R$.  The functor $\CechH_{I}^{0}$ is localization in the category of $R$-modules with respect
to the hereditary torsion theory of $I$-torsion modules.  
\end{corollary}

For this corollary to make sense, recall that a class of objects in an
abelian category $\cat{A}$ is a \textbf{hereditary torsion theory} if it
is closed under subobjects, extensions, quotient objects, and arbitrary
coproducts.  If $\cat{T}$ is a hereditary torsion theory, we define a
map $f$ to be a \textbf{$\cat{T}$-equivalence} if its kernel and
cokernel are in $\cat{T}$.  An object $X$ is called
\textbf{$\cat{T}$-local} if $\cat{A}(f,X)$ is an isomorphism for all
$\cat{T}$-equivalences $f$.  A \textbf{$\cat{T}$-localization} of an
object $M$ is a $\cat{T}$-local object $X$ together with a
$\cat{T}$-equivalence $M\xrightarrow{}X$.  When $\cat{T}$-localizations
exist, they are unique up to unique isomorphism and are functorial.

\begin{proof}
 One can easily check that the class $\cat{T}$ of $I$-torsion modules
 is a hereditary torsion theory.  It is clear from parts~(a) and~(c)
 of Proposition~\ref{prop-Cech} that the map
 $M\xrightarrow{}\CechH_{I}^{0}(M)$ has $I$-torsion kernel and
 cokernel, so is a $\cat{T}$-equivalence.  It remains to show that
 $\CechH_{I}^{0}(M)$ is $\cat{T}$-local.  By factoring a
 $\cat{T}$-equivalence into an injection followed by a surjection, we
 see that this boils down to showing that $\CechH_{I}^{0}(M)$ has no
 $I$-torsion and that $\Ext^{1}_{R}(T,\Cech H_{I}^{0}(M))=0$ for all
 $I$-torsion modules $T$.  The first part is part~(c) of
 Proposition~\ref{prop-Cech}.  For the second part, suppose we have an
 extension
\[
0 \xrightarrow{} \CechH_{I}^{0}(M) \xrightarrow{} X \xrightarrow{} T
\xrightarrow{} 0
\]
where $T$ is $I$-torsion.  Applying the left exact idempotent functor
$\CechH_{I}^{0}$, we get an isomorphism $\CechH_{I}^{0}(M)\xrightarrow{}\CechH_{I}^{0}(X)$.  Thus the composite 
\[
X \xrightarrow{} \CechH_{I}^{0}(X) \cong \CechH_{I}^{0}(M)
\]
defines a splitting of our extension.  Thus $\CechH_{I}^{0}(M)$ is
$\cat{T}$-local as required.  
\end{proof}

We also need to know that $\CechH^{*}_{I}$ are the right derived
functors of $\CechH^{0}$.  This seems to require some hypotheses on
$I$.

\begin{theorem}\label{thm-local}
 Suppose $I$ is an ideal in a commutative ring $R$, generated by a
 regular sequence $(\alpha_{0},\dotsc,\alpha_{n})$, in which each
 element $\alpha_i$ is not a zero-divisor.  Then
 \[ H_{I}^{k}(M) = \colim_J \Ext_R^k(R/J,M), \]
 where $J$ runs over ideals $J\leq I$ with $\sqrt{J}=\sqrt{I}$.
 Moreover, $H_{I}^{k}$ is the $k$'th right derived functor of
 $H_{I}^{0}$ and $\CechH_{I}^{k}$ is the $k$'th right derived functor
 of $\CechH_{I}^{0}$.
\end{theorem}

\begin{proof}[Proof of Theorem~\ref{thm-local}]
 Put $I_r=(\alpha_{0}^r,\dotsc,\alpha_{n}^r)$; these ideals are
 evidently cofinal among the $J$'s.  Let $K^\bt_r$ be the usual
 (unstable) Koszul complex for $I_r$, which is the tensor product over
 $j$ of the complexes $(R\xrightarrow{\alpha_j^r}R)$.  As our sequence of
 generators is regular, this is a finite resolution of $R/I_r$ by
 finitely generated free modules.  Now let $DK^\bt_r$ be the dual of
 $K^\bt_r$, which is naturally thought of as the tensor product of the
 complexes $R\xrightarrow{}R.\alpha_j^{-r}$.  (In fact $DK^\bt_r$ is
 isomorphic to $K^\bt_r$, up to a degree shift in the graded case.)
 It is clear that the stable Koszul complex $K^\bt(I)$ is the colimit
 of the complexes $DK^\bt_r$, so 
 \begin{align*}
  H^*_I(M) &= \colim_r H^*(DK^\bt_r\otimes_R M) \\
           &= \colim_r H^*\Hom_R(K^\bt_r,M) \\
           &= \colim_r \Ext^*_R(R/I_r,M) \\
           &= \colim_J \Ext^*_R(R/J,M).
 \end{align*}
 It is immediate from this that $H^i_I(M)=0$ if $i>0$ and $M$ is
 injective.  Using part~(a) of Proposition~\ref{prop-Cech}, we see
 that $\CechH^i_I(M)=0$ as well.  It now follows formally from the
 long exact sequences in Proposition~\ref{prop-Cech} that $H_{I}^{k}$
 is the $k$'th right derived functor of $H_{I}^{0}$, and
 $\CechH_{I}^{k}$ is the $k$'th right derived functor of
 $\CechH_{I}^{0}$.
\end{proof}

We can now investigate the functors $H_{I}^{*}$ and $\CechH_{I}^{*}$
restricted to the category of comodules, proving part~(1) of
Theorem~\ref{main-B}.

\begin{theorem}\label{thm-local-cohomology}
Suppose $M$ is a $BP_{*}BP$-comodule.  Then there are natural
isomorphisms $T_{n}^{k}(M)\cong H_{I_{n+1}}^{k}(M)$ and $L_{n}^{k}M\cong
\CechH_{I_{n+1}}^{k}(M)$.  
\end{theorem}

\begin{proof}
We first show that $H_{I_{n}}^{k}(M)=0$ for all injective
$BP_{*}BP$-comodules $M$ and $k>0$.  This does not follow from
Theorem~\ref{thm-local} because injective comodules need not be
injective as $BP_{*}$-modules.   We proceed by induction on
$n$, using the spectral sequence
\[
H_{v_{n}}^{s}H_{I_{n}}^{t}(M)\Rightarrow H_{I_{n+1}}^{s+t}(M).  
\]
discussed in~\cite[Section~2]{greenlees-may-completions}.
By induction, the $E_{2}$-term of this spectral sequence is
$H_{v_{n}}^{s}H_{I_{n}}^{0}(M)$.  In degree $s=1$, this is 
\[ v_{n}^{-1}H_{I_{n}}^{0}(M)/H_{I_{n}}^{0}(M). \]
But Theorem~\ref{thm-injective-torsion} shows that
$H_{I_{n}}^{0}M=T_{n-1}M$ is an injective $BP_{*}BP$-comodule, which of
course is $v_{n-1}$-torsion.  Proposition~\ref{prop-injective-divisible}
then shows that $H_{I_{n}}^{0}M$ is $v_{n}$-divisible.  Hence the
$E_{2}$-term of our spectral sequence is
$H_{v_{n}}^{0}H_{I_{n}}^{0}(M)=H_{I_{n+1}}^{0}(M)$ concentrated in
bidegree $(0,0)$, completing the proof. 

Now suppose $M$ is an arbitrary $BP_{*}BP$-comodule.  Take an resolution
$I_{*}$ of $M$ by injective $BP_{*}BP$-comodules.  By definition,
$T_{n}^{k}(M)\cong H_{-k}(T_{n}I_{*})$.  On the other hand, applying
$H_{I_{n+1}}^{*}$, which we have just seen vanishes on injective
comodules, shows that 
\[
H_{I_{n+1}}^{k}(M)\cong H_{-k}(H_{I_{n+1}}^{0}I_{*}) \cong H_{-k}(T_{n}I_{*})
\]
as well.  

Similarly, $L_{n}^{k}(M)\cong H_{-k}(L_{n}I_{*})$, which is isomorphic
to $H_{-k}(I_{*}/T_{n}I_{*})$ by Corollary~\ref{cor-injective-Ln}.  Now
suppose $N$ is an injective $BP_{*}BP$-comodule.  The exact sequence
\[
0 \xrightarrow{} H_{I_{n+1}}^{0}(N)
  \xrightarrow{} N
  \xrightarrow{} \CechH_{I_{n+1}}^{0}(N)
  \xrightarrow{} H_{I_{n+1}}^{1}(N)
  \xrightarrow{} 0
\]
of Proposition~\ref{prop-Cech} together with the fact that
$H_{I_{n+1}}^{k}(N)=0$ for $k>0$ implies that
$\CechH_{I_{n+1}}^{0}(N)\cong N/T_{n}N$.  Also,
\[
\CechH_{I_{n+1}}^{k}(N) \cong H_{I_{n+1}}^{k+1}(N)=0
\]
for $k>0$.  Hence, applying $\CechH_{I_{n+1}}^{*}$ to $I_{*}$, we find
that $\CechH_{I_{n+1}}^{k}(M)\cong H_{-k}(I_{*}/T_{n}I_{*})$,
completing the proof.  
\end{proof}

We can now give the promised converse to Lemma~\ref{lem-local-no-BP}.  

\begin{corollary}\label{cor-local-Cech}
A $BP_{*}BP$-comodule $M$ is $L_{n}$-local if and only if 
\[
\Hom_{BP_{*}}^{*}(BP_{*}/I_{n+1},M)
=\Ext^{1,*}_{BP_{*}}(BP_{*}/I_{n+1},M)=0.  
\]
\end{corollary}

\begin{proof}
The if direction is Lemma~\ref{lem-local-no-BP}.  For the only if
direction, suppose $M$ is $L_{n}$-local.  Then $M$ is also local with
respect to the hereditary torsion theory of $I_{n+1}$-torsion
$BP_{*}$-modules, in view of Theorem~\ref{thm-local-cohomology} and
Corollary~\ref{cor-Cech-local}.  This means that, for any
$I_{n+1}$-torsion module $T$, we have 
\[
\Hom_{BP_{*}}(T,M) = \Ext^{1}_{BP_{*}}(T,M)=0.
\]
Applying this to $BP_{*}/I_{n+1}$ and all its suspensions gives the
desired result.  
\end{proof}

\section{The spectral sequence}\label{sec-spectral}

The object of this section is to prove Theorem~\ref{main-A}.  That is,
we construct a spectral sequence converging to $BP_{*}L_{n}X$ whose
$E_{2}$-term consists of the derived functors $L^{s}_{n}(BP_{*}X)$.
Analogously, let $C_{n}X$ denote the fiber of $X\xrightarrow{}L_{n}X$.
We construct a spectral sequence converging to $BP_{*}C_{n}X$ whose
$E_{2}$-term consists of the derived functors $T^{s}_{n}(BP_{*}X)$.  Our
method is based on Devinatz' construction of the modified Adams spectral
sequence in~\cite[Section~1]{devinatz-morava-module}.  

\begin{definition}\label{defn-injective}
Define a functor $D$ from injective $BP_{*}BP$-comodules to (the
homotopy category of) spectra as
follows.  Given an injective $BP_{*}BP$-comodule $I$, consider the
functor $D_{I}$ from spectra to abelian groups defined by 
\[
D_{I}(X) = \Hom_{BP_{*}BP}(BP_{*}X, I).  
\]
Then $D_{I}$ is a cohomology functor, so there is a unique spectrum
$D(I)$ such that there is a natural isomorphism
\[
D_{I}(X)\cong [X, D(I)].
\]
\end{definition}

The reason for the letter $D$ is that $D_{I}$ is a sort of duality
functor, built along the lines of Brown-Comenetz
duality~\cite{brown-comenetz}.  Also note that we are considering
cohomology functors as exact functors to abelian groups; we recover the
usual graded cohomology functor by $D_{I}^{t}(X)=D_{I}(\Sigma^{t}X)$.

The following theorem is a special case of Theorem~1.5
of~\cite{devinatz-morava-module}.  

\begin{theorem}\label{thm-spectral-BP}
Suppose $I$ is an injective $BP_{*}BP$-comodule.  Then there is a
natural isomorphism $BP_{*}D(I)\cong I$.  
\end{theorem}

This isomorphism of course corresponds to the identity map of $D(I)$
under the isomorphism 
\[
[D(I), D(I)] \cong \Hom_{BP_{*}BP}(BP_{*}D(I),I).
\]

We need to know how the $D(I)$ behave under localization.  

\begin{proposition}\label{prop-spectral-Ln}
Suppose $I$ is an injective $BP_{*}BP$-comodule.  Then the natural map
$I\xrightarrow{}L_{n}I$ induces an isomorphism 
\[
L_{n}D(I) \xrightarrow{} D(L_{n}I).  
\]
\end{proposition}

\begin{proof}
Recall that $L_{n}I$ is again injective, by
Corollary~\ref{cor-injective-Ln}.  We first note that $D(L_{n}I)$ is
$E(n)$-local.  Indeed, if $E(n)_{*}(X)=0$, then $BP_{*}(X)$ is all
$v_{n}$-torsion.  Since $L_{n}I$ has no $v_{n}$-torsion, we have
\[
[X, D(L_{n}I)]\cong \Hom_{BP_{*}BP}(BP_{*}X, L_{n}I)=0. 
\] 
Thus $D(L_{n}I)$ is indeed $L_{n}$-local.  

On the other hand, the map $D(I)\xrightarrow{}D(L_{n}I)$ induces the map
$I\xrightarrow{}L_{n}I$ on $BP_{*}$-homology, by
Theorem~\ref{thm-spectral-BP}.  Since $L_{n}I\cong I/T_{n}I$ by
Corollary~\ref{cor-injective-Ln}, this map becomes an isomorphism after
applying $\Phi_{*}$, and so $D(I)\xrightarrow{}D(L_{n}I)$ is an
$E(n)$-equivalence.  
\end{proof}

\begin{corollary}\label{cor-spectral-Cn}
Suppose $I$ is an injective $BP_{*}BP$-comodule.  Then the natural map
$T_{n}I\xrightarrow{}I$ induces an isomorphism 
\[
D(T_{n}I) \xrightarrow{} C_{n}D(I).  
\]
\end{corollary}

\begin{proof}
Note that $D(T_{n}I)$ makes sense since $T_{n}I$ is an injective
comodule by Theorem~\ref{thm-injective-torsion}.  Since
$BP_{*}(D(T_{n}I))\cong T_{n}I$, one easily sees that $D(T_{n}I)$ is
$E(n)$-acyclic.  Therefore, the map $D(T_{n}I)\xrightarrow{}D(I)$
induced by the inclusion $T_{n}I\xrightarrow{}I$ induces the desired map
\[
D(T_{n}I) \xrightarrow{} C_{n}D(I).  
\]
This map is an isomorphism on $BP_{*}(-)$ by
Proposition~\ref{prop-spectral-Ln}, and one can check that both sides
are $BP$-local, so it is an isomorphism.
\end{proof}

We can now build our spectral sequences, following the standard approach
used by Ravenel in~\cite[Section~2.1]{ravenel}.  Suppose $X$ is a
spectrum, and let $C=BP_{*}X$.  Choose an injective resolution
\[
0 \xrightarrow{} C \xrightarrow{\eta } I_{0} \xrightarrow{\tau _{0}} I_{1}
\xrightarrow{\tau _{1}} \dotsb 
\]
of $C$ in the category of $BP_{*}BP$-comodules.  Let $\eta _{s}\mathcolon
C_{s}\xrightarrow{}I_{s}$ denote the kernel of $\tau _{s}$, so that
$\eta _{0}=\eta$.

The following lemma is easily proved by induction on $s$, and is
implicit in~\cite[Section~1]{devinatz-morava-module}.  

\begin{lemma}\label{lem-spectral-construction}
Let $X$ be a spectrum and choose an injective resolution of $BP_{*}X$ as
above.  Then there is a tower 
\[
\begin{CD}
X=X_{0} @<g_{0}<< X_{1} @<g_{1}<< X_{2} @<g_{2}<< \dotsb \\
@. @VVf_{0}V @VVf_{1}V \\
\ @. K_{0} @. K_{1} 
\end{CD}
\]
over $X$ satisfying the following properties. 
\begin{enumerate}
\item [(a)] $K_{s}=\Sigma^{-s}D(I_{s})$.  
\item [(b)] $X_{s+1}$ is the fiber of $f_{s}$.  
\item [(c)] $BP_{*}X_{s}\cong \Sigma^{-s}C_{s}$.  
\item [(d)] The map $f_{s}$ is induced by the inclusion
$C_{s}\xrightarrow{}I_{s}$.  
\item [(e)] $BP_{*}g_{s}=0$, and the boundary map
$K_{s}\xrightarrow{}\Sigma X_{s+1}$ induces the surjection
$\Sigma^{-s}I_{s}\xrightarrow{}\Sigma^{-s}C_{s+1}$ on
$BP_{*}$-homology.  
\end{enumerate}
\end{lemma}

We can now construct our spectral sequences.  The following theorem
is Theorem~\ref{main-A} except for the statements about convergence. 

\begin{theorem}\label{thm-spectral-build-Ln}
Let $X$ be a spectrum.  There is a natural spectral sequence
$E_{*}^{**}(X)$ with $d_{r}\mathcolon
E_{r}^{s,t}\xrightarrow{}E_{r}^{s+t,t+r-1}$ and $E_{2}$-term
$E_{2}^{s,t}(X)\cong (L^{s}_{n}BP_{*}X)_{t}$.  This is a spectral
sequence of $BP_{*}BP$-comodules, in the sense that $E_{r}^{s,*}$ is a
graded $BP_{*}BP$-comodule for all $r\geq 2$ and $d_{r}\mathcolon
E_{r}^{s,*}\xrightarrow{}E_{r}^{s+r,*}$ is a $BP_{*}BP$-comodule map of
degree $r-1$.  Furthermore, every element in $E_{2}^{0,*}$ that comes
from $BP_{*}X$ is a permanent cycle.
\end{theorem}

\begin{proof}
Begin with the tower of Lemma~\ref{lem-spectral-construction} and apply
$L_{n}$.  We get the tower below.
\begin{equation}\label{eq-ln-tower}
\begin{CD}
L_{n}X=L_{n}X_{0} @< L_{n}g_{0}<<  L_{n}X_{1} @< L_{n}g_{1}<<
L_{n}X_{2} @< L_{n}g_{2}<< \dotsb \\ 
@. @VV L_{n}f_{0}V @VV L_{n}f_{1}V  \\
\ @.  L_{n}K_{0} @.  L_{n}K_{1} 
\end{CD}
\end{equation}
By applying $BP_{*}$-homology, we get an associated exact couple and
spectral sequence.  That is, we let $D_{1}^{s,t}=BP_{t-s}L_{n}X_{s}$ and
$E_{1}^{s,t}=BP_{t-s}L_{n}K_{s}.$  We take 
\[
i_{1}=BP_{t-s}L_{n}g_{s} \mathcolon
D_{1}^{s+1,t+1}\xrightarrow{}D_{1}^{s,t} \text{ and }
j_{1}=BP_{t-s}L_{n}f_{s}\mathcolon D_{1}^{s,t} \xrightarrow{}
E_{1}^{s,t}
\]
and we take
\[
k_{1}\mathcolon E_{1}^{s,t} \xrightarrow{} D_{1}^{s+1,t}
\]
to be $BP_{t-s}$ of the boundary map $L_{n}K_{s}\xrightarrow{}\Sigma
L_{n}X_{s+1}$.  

Note that this is an exact couple in the category of
$BP_{*}BP$-comodules, in that each $D_{1}^{s,*}$ and $E_{1}^{s,*}$ is a
graded $BP_{*}BP$-comodule and the maps $i_{1},j_{1},k_{1}$ are maps of
comodules.  It follows that the spectral sequence is a spectral sequence
of $BP_{*}BP$-comodules.

Now, by combining Theorem~\ref{thm-spectral-BP} and
Proposition~\ref{prop-spectral-Ln}, we find 
\[
E_{1}^{s,t} \cong BP_{t-s}(L_{n}\Sigma^{-s}D(I_{s})) \cong
BP_{t}D(L_{n}I_{s}) \cong (L_{n}I_{s})_{t}.  
\]
To compute the first differential $d_{1}$, note that we have the
commutative diagram below.  
\[
\begin{CD}
K_{s} @>>> \Sigma X_{s+1} @>>> \Sigma K_{s+1} \\
@VVV @VVV @VVV \\
L_{n}K_{s} @>>> L_{n}(\Sigma X_{s+1}) @>>> L_{n}(\Sigma K_{s+1})
\end{CD}
\]
The map on $BP_{t-s}$-homology induced by the bottom composite is
$d_{1}$.  The map on $BP_{t-s}$-homology induced by the top composite is
$\tau _{s}$, by Lemma~\ref{lem-spectral-construction}.  The
outside vertical maps are surjective in $BP_{*}$-homology, by
Proposition~\ref{prop-spectral-Ln} and
Corollary~\ref{cor-injective-Ln}.  It follows that
$d_{1}=L_{n}\tau_{s}$.  Therefore, the $E_{2}$-term of our
spectral sequence is 
\[
E_{2}^{s,t} \cong H^{s}(L_{n}I_{*})_{t} \cong (L^{s}_{n}BP_{*}X)_{t},
\]
as required.  

The naturality of the spectral sequence follows in the usual way.  That
is, a map of spectra $X\xrightarrow{}Y$ induces a map
$BP_{*}X\xrightarrow{}BP_{*}Y$.  This can be lifted, nonuniquely, to a
map of injective resolutions and so to a map of the towers of
Lemma~\ref{lem-spectral-construction}.  This map induces a map of
spectral sequences which is the evident map 
\[
L^{s}_{n}(BP_{*}X) \xrightarrow{} L^{s}_{n}(BP_{*}Y)
\]
on the $E_{2}$-terms.  This map is independent of the choice of map of
injective resolutions, and so is functorial.  This also shows that our
spectral sequence is independent of the choice of injective resolution
(from $E_{2}$ on).  

To complete the proof, we must show that every element in $E_{2}^{0,*}$
that comes from $BP_{*}X$ is a permanent cycle.  To see this, note that
there is a map from the tower of Lemma~\ref{lem-spectral-construction}
to the tower~\ref{eq-ln-tower} induced by $L_{n}$.  Applying
$BP_{*}$-homology to the tower of Lemma~\ref{lem-spectral-construction}
gives us a spectral sequence with $E_{2}^{s,t}=0$ if $s>0$ and
$E_{2}^{0,t}=BP_{t}X$.  The map from this spectral sequence to our
spectral sequence immediately gives the desired result.  
\end{proof}

We have an analogous theorem for $C_{n}$. 

\begin{theorem}\label{thm-spectral-build-Cn}
Let $X$ be a spectrum.  There is a natural spectral sequence
$E_{*}^{**}(X)$ with $d_{r}\mathcolon
E_{r}^{s,t}\xrightarrow{}E_{r}^{s+t,t+r-1}$ and $E_{2}$-term
$E_{2}^{s,t}(X)\cong (T^{s}_{n}BP_{*}X)_{t}$.  This is a spectral
sequence of $BP_{*}BP$-comodules, in the sense that $E_{r}^{s,*}$ is a
graded $BP_{*}BP$-comodule for all $r\geq 2$ and $d_{r}\mathcolon
E_{r}^{s,*}\xrightarrow{}E_{r}^{s+r,*}$ is a $BP_{*}BP$-comodule map of
degree $r-1$.  
\end{theorem}

\begin{proof}
Begin with the tower of Lemma~\ref{lem-spectral-construction} and apply
$C_{n}$ to get the tower below.  
\begin{equation}\label{eq-cn-tower}
\begin{CD}
C_{n}X=C_{n}X_{0} @< C_{n}g_{0}<<  C_{n}X_{1} @< C_{n}g_{1}<<
C_{n}X_{2} @< C_{n}g_{2}<< \dotsb \\ 
@. @VV C_{n}f_{0}V @VV C_{n}f_{1}V \\
\ @.  C_{n}K_{0} @.  C_{n}K_{1} 
\end{CD}
\end{equation}
Apply $BP_{*}$-homology to get an associated exact couple and spectral
sequence, as in the proof of Theorem~\ref{thm-spectral-build-Ln}.  This
time the $E_{1}$ term will be 
\[
E_{1}^{s,t}\cong BP_{t-s}C_{n}K_{s} \cong (T_{n}I_{s})_{t},
\]
using Corollary~\ref{cor-spectral-Cn}.  The identification of the
$E_{2}$-term uses the commutative diagram below.  
\[
\begin{CD}
C_{n}K_{s} @>>> C_{n}(\Sigma X_{s+1}) @>>> C_{n}(\Sigma K_{s+1}) \\
@VVV @VVV @VVV \\
K_{s} @>>> \Sigma X_{s+1} @>>> \Sigma K_{s+1} \\
\end{CD}
\]
The vertical maps are injective on $BP_{*}$-homology by
Corollary~\ref{cor-spectral-Cn} and
Theorem~\ref{thm-injective-torsion}.  Thus $d_{1}$, which is the effect
on $BP_{t-s}$-homology of the top horizontal composite, is
$T_{n}\tau_{s}$.  Hence we get the desired $E_{2}$-term and naturality,
as in Theorem~\ref{thm-spectral-build-Ln}.  
\end{proof}

We must now prove that our spectral sequences converge, strongly and
conditionally.  This essentially boils down to showing that the homotopy
inverse limits of the towers~\ref{eq-ln-tower} and~\ref{eq-cn-tower} are
trivial.  The plan of the proof is very simple; in the original tower of
Lemma~\ref{lem-spectral-construction}, we have $BP_{*}g_{s}=0$.  Hence
$E(n)_{*}(L_{n}g_{s})=E(n)_{*}g_{s}=0$ as well by Landweber exactness.
Now we just apply the following theorem.  

\begin{theorem}\label{thm-spectral-composition}
Given $n\geq 0$, there exists an $N$ such that every composite
\[
g=f_{N}\circ f_{N-1}\circ \dotsb \circ f_{1}
\]
of maps of spectra such that $E(n)_{*}f_{i}=0$ for all $i$ has
$L_{n}g=0$.
\end{theorem}

This theorem was certainly known to Hopkins and probably others.

\begin{proof}
Use the modified Adams spectral sequence
\[
E_{2}^{s,t}=\Ext^{s,t}_{E(n)_{*}E(n)}(E(n)_{*}X, E(n)_{*}Y) \Rightarrow
[X, L_{n}Y]_{t-s}
\]
of Devinatz~\cite{devinatz-morava-module}.  It was proved
in~\cite[Proposition~6.5]{hovey-strickland} that there are integers
$r_{0}$ and $s_{0}$, independent of $X$ and $Y$, such that
$E_{r}^{s,t}=0$ whenever $r\geq r_{0}$ and $s\geq s_{0}$.  Take
$N=s_{0}$.  Then the composite $g$ is represented by an element in
$E_{2}^{s,t}$ for some $s\geq N$.  Therefore $g$ must be represented by
some element in $E_{\infty}^{s,t}$ for some $s\geq N$, so $g=0$.
\end{proof}

The following corollary completes the proof of Theorem~\ref{main-A}.  

\begin{corollary}\label{cor-spectral-convergence-Ln}
The spectral sequence of Theorem~\ref{thm-spectral-build-Ln} converges
strongly and conditionally to $BP_{*}L_{n}X$.  
\end{corollary}

\begin{proof}
In view of Theorem~\ref{thm-spectral-composition}, the composites
$L_{n}X_{k+s}\xrightarrow{}L_{n}X_{k}$ in the tower~\ref{eq-ln-tower}
are trivial for large $s$.  Hence $\lim_{s}
BP_{*}L_{n}X_{s}=\lim^{1}_{s} BP_{*}L_{n}X_{s} =0,$ and so the spectral
sequence converges conditionally to $BP_{*}L_{n}X$~\cite{boardman}.  On
the other hand, it is clear that $\lim^{1}_{r}E_{r}^{s,t}=0$, since we
have a horizontal vanishing line.  Thus, the spectral sequence converges
strongly as well~\cite[Theorem~7.3]{boardman}.
\end{proof}

We also want to know that the other spectral sequence we have
constructed converges. 

\begin{corollary}\label{cor-spectral-convergence-Cn}
The spectral sequence of Theorem~\ref{thm-spectral-build-Cn} converges
strongly and conditionally to $BP_{*}C_{n}X$.  
\end{corollary}

\begin{proof}
We have a cofiber sequence
$C_{n}X_{s}\xrightarrow{}X_{s}\xrightarrow{}L_{n}X_{s}$ of towers, where
$X_{s}$ denotes the tower of Lemma~\ref{lem-spectral-construction}.  By
applying $BP_{*}$, we get an exact sequence of towers 
\[
BP_{*+1}L_{n}X_{s} \xrightarrow{} BP_{*}C_{n}X_{s} \xrightarrow{}
BP_{*}X_{s} \xrightarrow{} BP_{*}L_{n}X_{s}.  
\]
We have just seen, in Corollary~\ref{cor-spectral-convergence-Ln}, that
the towers $BP_{*+1}L_{n}X_{s}$ and $BP_{*}L_{n}X_{s}$ are pro-trivial.
It follows that the tower $BP_{*}C_{n}X_{s}$ is pro-isomorphic to the
tower $BP_{*}X_{s}$.  But the tower $BP_{*}X_{s}$ is obviously
pro-trivial by Lemma~\ref{lem-spectral-construction}, so the tower
$BP_{*}C_{n}X_{s}$ is also pro-trivial.  Hence $\lim_{s}
BP_{*}C_{n}X_{s} \cong \lim^{1}_{s} BP_{*}C_{n}X_{s} =0,$ and so the
spectral sequence of Theorem~\ref{thm-spectral-build-Cn} converges
conditionally.  Since it has a horizontal vanishing line, it also
converges strongly~\cite[Theorem~7.3]{boardman}.
\end{proof}

We close the paper by considering the spectral sequence of
Theorem~\ref{main-A} in case $X=S^{0}$ and $n>0$.  In that case, we have
$E_{2}^{0,*}\cong BP_{*}$ and $E_{2}^{n,*}=BP_{*}/I_{n+1}^{\infty}$, by
Corollary~\ref{cor-derived-standard}.  The only possible differential is
$d_{n}$, but this must be trivial since $E_{2}^{0,*}$ must consist of
permanent cycles by Theorem~\ref{main-A}.  Thus our spectral sequence
degenerates to the short exact sequence of comodules
\[
0 \xrightarrow{} \Sigma^{-n} BP_{*}/I_{n+1}^{\infty} \xrightarrow{} BP_{*}L_{n}S^{0}
\xrightarrow{} BP_{*} \xrightarrow{} 0.
\]
A splitting of this sequence is given by the map
$BP_{*}\xrightarrow{}BP_{*}L_{n}S^{0}$ induced by
$S^{0}\xrightarrow{}L_{n}S^{0}$.  Hence we recover Ravenel's computation
of $BP_{*}L_{n}S^{0}$~\cite[Theorem~6.2]{rav-loc}.


\providecommand{\bysame}{\leavevmode\hbox to3em{\hrulefill}\thinspace}
\providecommand{\MR}{\relax\ifhmode\unskip\space\fi MR }
\providecommand{\MRhref}[2]{%
  \href{http://www.ams.org/mathscinet-getitem?mr=#1}{#2}
}
\providecommand{\href}[2]{#2}

\end{document}